\documentclass[12pt]{article}

\def\chi{{\mathcal{X}}}

\def\cald{{\mathcal{D}}}
\def\calx{{\mathcal{X}}}
\def\call{{\mathcal{L}}}

\def\({\left(}
\def\){\right)}
\def\pf{\n{\bf Proof.} }
\def\vsp{\vspace*{1,5mm}\\ }

\def\bk{\bigskip }

\def\sk{\smallskip }

\def\n{\noindent }
\def\dd{\displaystyle}

\def\D{{\Delta}}

\def\barr{\begin{array}}
\def\earr{\end{array}}

\def\bit{\begin{itemize}}
\def\itemi{\item[{\rm(i)}]}
\def\itemii{\item[{\rm(ii)}]}

\def\eit{\end{itemize}}

\def\D{{\Delta}}

\usepackage{amsmath, amsthm, amscd, amsfonts, amssymb}

\newtheorem{theorem}{Theorem}[section]

\theoremstyle{definition}
\newtheorem{definition}[theorem]{Definition}

\newtheorem{remark}[theorem]{Remark}

\def\1{^{-1}}
\def\vsp{\vspace*{2mm}\\ }

\def\calf{{\mathcal{F}}}

\def\calx{{\mathcal{X}}}

\def\rr{{\mathbb{R}}}

\def\9{{\infty}}

\def\lbb{{\lambda}}

\def\b{{\beta}}

\def\g{{\gamma}}

\def\ov{\overline}
\def\vf{{\varphi}}

\def\ooo{{\Omega}}
\def\pp{{\partial}}
\def\D{{\Delta}}
\def\vp{{\varepsilon}}
\def\barr{\begin{array}}
\def\earr{\end{array}}

\def\dd{\displaystyle}
\def\bk{\bigskip }
\def\sk{\smallskip}
\def\n{\noindent }

\def\vsp{\vspace*{2mm}\\ }
\def\ff{\forall }
\def\({\left(}
\def\){\right)}
\def\<{\left<}
\def\>{\right>}
\def\one{{\bf1}}

\title{Probabilistic representation for solutions to~nonlinear  Fokker--Planck equations}
\author{Viorel Barbu\thanks{Octav Mayer Institute of
Mathematics of Romanian  Academy, Ia\c si, Romania} \and
Michael R\"ockner\thanks{Fakult\"at f\"ur Mathematik, Universit\"at Bielefeld, 33615, Bielefeld, Germany}}
\date{}

\begin{document}
\maketitle
\begin{abstract}
\n One obtains a probabilistic representation for the entropic ge\-ne\-ra\-lized solutions to a nonlinear  Fokker--Planck equation in $\rr^d$ with multi\-valued nonlinear diffusion term as density pro\-ba\-bi\-li\-ties of solutions to a nonlinear stochastic differential equation.  The case of a non\-linear Fokker--Planck equation with linear space dependent drift  is also \mbox{studied.}\sk
\end{abstract}

\section{Introduction}\label{s1}
Consider here the nonlinear Fokker--Planck equation (NFPE)
\begin{equation}\label{e1.1}
\barr{l}
u_t(t,x)+{\rm div}(b(u)u)-\D \b(u)=0\ \mbox{ in }(0,\9)\times\rr^d,\\
u(0,x)=u_0(x),\ \ x\in\rr^d,\ d\ge1,\earr\end{equation}under the following assumptions
\bit\itemi $\b:\rr\to2^\rr$ is  a maximal monotone (multivalued) function  $\b(0)=0.$
\itemii $b\in C_b(\rr;\rr^d)$.
\eit
Assumption (i) means that $(\eta_1-\eta_2)(u_1-u_2)\ge0$ for all $(u_i,\eta_i)\in\rr\times\rr,$ such that $\eta_i\in\b(u_i)$, $i=1,2,$ and the range of the mapping $r\to r+\b(r)$ is all of $\rr$. In particular, this holds if $\b$ is continuous and monotonically nondecreasing. If $\b$ is a discontinuous, single-valued, monotonically increasing function on $\rr$, by filling the jumps $r_j$, that is, redefining $\b$ as $\b(r_j)=\b(r_j+0)-\b(r_j-0)$, one gets a maximal monotone $\b:\rr\to2^\rr.$

In the mean field theory and statistical mechanics,  equation
\eqref{e1.1} describes the
particle transport dynamics in disorded media and $u(t)$ is the probability density.

In general, a Fokker-Planck equation of the form \eqref{e1.1} is associated with a certain entropic functional. For instance, in $1{-}D$,    equation \eqref{e1.1} is derived from the entropic functional
$$S[u]=\int_\rr G(u(x))dx,$$
where $G\in C^2(0,\9),\ G'(0)=\9,\ G''(u)<0,$ $\ff u>0.$ Then the NFPE associated with processes with entropy $S$ is (see, e.g., \cite{7}, \cite{8})
$$P_t+(b(P)P-\g(G(P)-PG'(P))_{x})_x=0,$$where $\g$ is some positive constant.
This is an equation of the form \eqref{e1.1}, where $\b(u)\equiv G(u)-uG'(u)$. A notorious example is $b=1$ and $\beta(u)=$\break $\g\ln(1+u)$  (equation of classical bosons).

Another important example is $\b(t)\equiv rH(r-\rho_c)$, where $H$ is the Heaviside function and $r_\vp>0$. In this case, equation \eqref{e1.1} describes the self-organized criticality with the entropic functional (see \cite{3ab})
$$g(r)=\left\{\barr{ll}
\dd\frac r{\rho_c}\ln\(\frac r{\rho_c}\)&\mbox{ for }r>\rho_c,\vsp
0&\mbox{ for }r<\rho_c.\earr\right.$$

In general,   NFPE \eqref{e1.1} has not a classical solution, but under assumptions (i)--(ii) it has an {\it entropic generalized solutions} in the sense of S. Kru\v zkov (see~\cite{2}). (See Section 3.)

The main result of this paper, Theorem \ref{t3.1}, amounts to saying that, for   the single valued $\beta$ and $b$ as above,
the generalized entropic solution to NFPE \eqref{e1.1} can be represented as the law density  of a stochastic process $ Y$ which is (in the probabilistic sense) a weak solution to    a stochastic differential equation with the drift $b$ and the diffusion term  $\Phi(u)=(2\beta(u)/u\one_{[0,\9)}(u))^{1/2},$ that is,
\begin{equation}\label{11a}
dY(t)=b(u(t,Y(t))dt+\Phi(u(t,Y(t)))dW(t),\end{equation}
with a suitable modification for multivalued functions $\b$. Equation \eqref{11a} can be viewed as a stochastic model of NFPE \eqref{e1.1}. Pre\-viously, such a $1{-}D$ result was obtained in \cite{7a}, \cite{4}, \cite{5} for $b\equiv0$, $\b(u)=u^m$ and $\b(u)=H(u-u_c)u$, respectively. In the special case $\b(u)\equiv u$, one obtains the re\-pre\-sen\-ta\-tion of solutions to the linear Fokker-Planck equation as probability density of diffusion processes with the Boltzmann-Gibbs entropy $$S[u]=-\int_\rr u\log udx.$$
In \cite{13a}, the existence of a solution to a McKean type SDE of the form \eqref{11a} is studied and it is proved that its law density is a distributional solution to a corresponding Fokker-Planck equation.

The principal problem encountered with the probabilistic representation of entropic solutions to \eqref{e1.1} is their weak regularity. (See Definition \ref{d2.1} below.) To circumventing this, we developed here a new approach which, in  a few words, can be described as follows; one approximates \eqref{e1.1} by  a family of smooth parabolic problems for which one has such a representation and gets the result by passing to the limit in the corresponding stochastic equation.

The structure of the paper is as follows. In Section 2, we briefly sketch the standard technique to obtain probabilistic representations for solutions of    nonlinear Fokker-Planck equations. In Section 3, we recall the notion of entropic solutions and present a proof, via approximation for existence of a solution. In Section 4, we derive a probabilistic representation for the latter. In Section 5, we study the case when $b$ is independent of $u$, but depending on the spatial variable. In Section 6, an application to $1{-}D$ degenerate parabolic equation is given.

\bk\n{\bf Notation.} For $1\le p\le\9$, we denote by $L^p$    the space of  Lebesgue $p$ integrable  functions $u:\rr^d\to\rr$ with   norm denoted by  $|\cdot|_p$,   by  $H^1$  we denote the Sobolev space,  $H^{1}(\rr^d)=\left\{u\in L^2;\frac{\pp u}{\pp x_i}in L^2(\rr^d)\right\}$, and by $W^{1,p}_{\rm loc}$,   $H^1_{\rm loc}$,  the corresponding local spaces; $H^{-1}$ is the dual of $H^1$.   $C^k(\rr^d)$, $k=0,1,2$, is the space of $k$-times continuously differentiable functions. We denote by ${\rm Lip}(\rr)$ the space of Lipschitz functions on $\rr$ and by ${\rm Lip}_p^{\rm loc}(\rr)$ the space of locally Lipschitz functions on $\rr$. By $C^k(\rr^d;\rr^m)$, we denote the corresponding $\rr^m$-valued function  space. We set $C(\rr^d;\rr^m)=C^0(\rr^d;\rr^m)$ and let $C_b(\rr^d;\rr^m)$ denote the set of all bounded functions in $C(\rr^d;\rr^m)$. Furthermore, $C^k_0(\rr^d;\rr^m)$ denotes the space  of all the functions in $C^k(\rr^d;\rr^m)$ with compact support. We shall use also the notations: $x=(x_1,...,x_d)$, $b=(b_1,b_2,...,b_d)$,    ${\rm sign}_,r=\frac r{|r|}$ for $r\ne0$, ${\rm sign}\,0=0$, ${\rm div}\,a(u)=\sum\limits^d_{i=1}\frac\pp{\pp x_i}\,a_i(u)$, $\nabla u=\left\{\frac{\pp u}{\pp x_i}\right\}^d_{i=1}$, $\D u=\sum\limits^d_{i=1}\frac{\pp^2u}{\pp x^2_i},$ $u_t=\frac{\pp u}{\pp t}$. By $|\cdot|_d$, we denote the Euclidean norm of $\rr^d$, $d\ge1$.
Given a Hilbert space $X$,   $C([0,T];X)$ denote  the space of $X$-valued continuous functions on $[0,T]$.

\section{Probabilistic representation of solutions\\  to NFPE}
\setcounter{equation}{0}
\setcounter{theorem}{0}

Let $T\in(0,\9)$ and $a^{ij}:[0,T]\times\rr\times\rr^d\to\rr$, $b^i:[0,T]\times\rr\times\rr^d\to\rr$ be Borel-measurable maps such that $(a^{ij}(t,r,x))_{1\le i,j\le d}$ is a symmetric nonnegative matrix for all $(t,r,x)\in[0,T]\times\rr\times\rr^d$. Consider the following nonlinear Fokker-Planck equation
\begin{equation}
\label{2.1z}\barr{l}
\pp_tu=\pp_i\pp_j(a^{ij}(u)u)-\pp_j(b^j(u)u),\vsp
u(0,\cdot)=u_0,\earr
\end{equation}
where $u_0$ is a probability density (with respect to the Lebesgue measure on~$\rr^d$) and we used Einstein's summation convention. Equation \eqref{2.1z} is to be understood in the weak sense, i.e.,
\begin{equation}
\label{2.2z}\barr{l}
\dd\int^T_0\int_{\rr^d}
[\pp_t \vf(t,x)
+a^{ij}(t,u(t,x),x)\pp_i\pp_j\vf(t,x)\\
\qquad\qquad
+b^i(t,u(t,x),x)\pp_i\vf(t,x)]u(t,x)dx\,dt\\
\qquad\qquad=-\dd\int_{\rr^d}\vf(0,x)u_0(x)dx,\earr
\end{equation}
for all $\vf$ of the form
$$\vf(t,x)=f(t)g(x),\ (t,x)\in\rr\times\rr^d,$$
with $f\in C ^1([0,T];\rr),$ $f(T)=0,$ $g\in C^2_0(\rr^d)$, where, in \eqref{2.2z}, it is assumed that
\begin{eqnarray}
\label{2.3z}\dd\int^T_0\int_{\rr^d}(|a^{ij}(u)|+|b^i(u)|)|u|dx\,dt<\9,\ 1\le i,j\le d,
\end{eqnarray}(see \cite{15}).

Assume that a solution $u:[0,T]\times\rr^d\to\rr$ exists such that
$$\int_{\rr^d} u(t,x)dx=1\mbox{ for all }t\in[0,T],$$and $t\mapsto u(t,x)dx$ is weakly continuous on $[0,T]$. Then, as a consequence of the so-called {\it superposition principle}, more precisely Theorem 2.5 in \cite{15} (see also \cite{9a}), there exists a probability measure $P_{u_0}$ on $C([0,T];\rr^d)$ equipped with its Borel $\sigma$-algebra and its natural filtration generated by the evaluation maps $\pi_t$, $t\in[0,T]$, defined by
$$\pi_t(w):=w(t).\ w\in C([0,T],\rr^d),$$solving the martingale problem in the sense of \cite{15} for the Kolmogorov operator
$$L=\pp_t+a^{ij}(u)\pp_i\pp_j+b^i(u)\pp_i$$with the initial distribution $P_{u_0}\circ\pi_0^{-1}=u_0 dx$ and, moreover, with marginals $P_{u_0}\circ\pi^{-1}_t=u(t,x)dx.$ Then, by a standard result (see, e.g.,  \cite{14prim}, Theorem 2.6), there exists a $d$-dimensional $(\calf_t)$-Brownian motion $W(t),$ $t\in[0,T]$, on a stochastic basis $(\ooo,\calf,(\calf_t)_{t\ge0},Q)$ and a continuous $(\calf_t)$-progressively measurable map $Y:[0,T]\times\ooo\to\rr^d$ satisfying the following distribution dependent stochastic differential equation (SDDDE) (see \cite{16})
\begin{equation}
\label{24z}
dY(t)=b(t,u(t,Y(t)),Y(t))dt+\sqrt{2}\,
\sigma(t,u(t,Y(t)),Y(t))dW(t)
\end{equation}
and the law
\begin{equation}
\label{2.4z} Q\circ Y^{-1}=P_{u_0},
\end{equation}where $\sigma=((a^{ij})_{1\le i,j\le d})^{\frac12}.$ In particular, we have, for the marginals,
$$Q\circ (Y(t))^{-1}=u(t,x)dx,\ t\in[0,T].$$
So, we have obtained a probabilistic representation of the solution $u$ of \eqref{2.1z} (in the sense of \eqref{2.2z}), i.e., $u(t,\cdot)$ is the law density of a stochastic process $Y$ which is a weak solution to \eqref{2.4z}.

It is much harder to prove that the solution to \eqref{24z} is unique in law, provided its initial distribution is $u_0dx$, which would of course be very de\-si\-rable. For this,  one has to prove that the solutions to the linear  Fokker-Planck equation
\begin{equation}
\label{2.5z}\barr{l}
\pp_tv=\pp_i\pp_j(a^{ij}(u)v)-\pp_j(b^j(u)v)\vsp
v(0,\cdot)=u_0\earr
\end{equation}in the sense of \eqref{2.2z} with $u(t,x)dx\,dt$ replaced by $v(t,x)dx\,dt$. This was, however, achieved in certain cases (see \cite{4}, \cite{5}  and also \cite{15a}).

As explained in the introduction of this paper, we look at generalized (=~entropic) solutions for a special case of \eqref{2.1z}. In this case, the above approach  applies directly,  because the entropic solution to \eqref{e1.1} in the sense of \eqref{e2.1} below is also a distribution solution and has the properties required above. However,   the idea   to approximate entropic solutions by solutions of more regular equations satisfying \eqref{2.2z} was necessary in order to get the necessary continuity properties.  This procedure will be implemented in the following two sections.

\section{Generalized entropic solutions to NFPE}\label{s2}
\setcounter{equation}{0}
\setcounter{theorem}{0}

Let $b$ as in (ii) above and set $a(r):=b(r)r,$ $\ff r\in\rr.$

\begin{definition}\label{d2.1}{\rm A function $u\in L^1((0,T)\times\rr^d)\cap C([0,T];L^1_{\rm loc})$ is said to be an entropic  solution to NFPE \eqref{e1.1} if there is $\eta\in L^2(0,T;H^1_{\rm loc})$ such that $\eta(t,x)\in\b(u(t,x))$, a.e. $(t,x)\in(0,T)\times\rr^d$ and, for all $k\in\rr$,
		\begin{equation}
		\label{e2.1}
		\barr{l}
		\dd\frac{\pp}{\pp t}\,|u-k|+{\rm div}_x({\rm sign}(u-k)(a(u)-a(k)))\\
		\quad\!\!\qquad\qquad+\,{\rm div}_x(\nabla_x\eta\,{\rm sign}(u-k))\le0\mbox{ in }\cald'((0,T)\times\rr^d),\vsp
		u(0,x)=u_0(x),\ x\in\rr^d.\earr
		\end{equation}
Equivalently, the initial condition holds and
\begin{equation}
\label{e2.1prim}
\barr{l}
\dd\int^T_0\int_{\rr^d}(|u(t,x)-k|\vf_t(t,x)+{\rm sign}(u(t,x)-k)\vsp
\qquad a(u(t,x)-a(k))\cdot\nabla_x\vf(t,x)-{\rm sign}(u(t,x)-k)\vsp
\qquad \nabla_x \eta(t,x) \cdot\nabla_x\vf(t,x))dt\ge0\earr
\end{equation}
for all $\vf\in C^\9_0((0,T)\times\rr^d),\ \vf\ge0$ and all $k\in\rr$.}
\end{definition}
(In the following, we shall simply write $\nabla$ instead  of $\nabla_x$.)

If $\b=0$, this is just Kruzkov's entropic solution to the conservation law equation (see \cite{6}, \cite{9}).

The function $\eta$ associated with $u$ in Definition \ref{d2.1} will be called an entropic co-solution to \eqref{e1.1} and will be denoted by $\eta^u$.

We note that, if the entropic solution $u$ is in $L^\9((0,T)\times\rr^d)$, then it is also a solution in sense of distributions to \eqref{e1.1}, that is,
$$\int^T_0\int_{\rr^d}(u\vf_t+a(u)\cdot\nabla
\vf
-\nabla_x\eta^u\cdot\nabla \vf)dt\,dx=0,\ \ff\vf\in C^\9_0((0,T)\times\rr^d).$$
This follows by taking in \eqref{e3.2} $k=|u|_\9+1$ and, respectively, $k=-|u|_\9-1.$

In \cite{2},    the existence of an entropic solution $u\in C([0,T];L^1)$ was proven via the Crandall and Liggett generation theorem in $L^1$. Here we shall give a direct constructive proof based on an appropriate smooth approximating equation. However, as seen in Remark \ref{r2.1} below, the construction  leads to the same concept of solution.

Namely, for $\vp\in(0,1]$, consider the equation
\begin{equation}\label{e2.11}
 \barr{l}
(u_\vp)_t+{\rm div}(a_\vp(u_\vp))-\D(\b_\vp(u_\vp)+ \vp u_\vp)=0 \mbox{ in }(0,\9){\times}\rr^d,\\
u_\vp(0,x)=u_0(x),\ x\in\rr^d.\earr \end{equation}
Here $a_\vp(u)\equiv b_\vp(u)u$ and the functions $\b_\vp ,b_\vp$ are defined as follows.
\begin{eqnarray}
\b_\vp(r)&\!\!\!=\!\!\!&
\frac1\vp\,(r-(1+\vp\b)^{-1}r)
\in \b((1+\vp\b)^{-1}r),\, \ff r\in\rr,\ \ \ \ \
\label{e3.4z}\\[2mm]
b_\vp(r)&\!\!\!=\!\!\!&\frac1\vp\int_{\rr}\frac{b(\theta)}{1+\vp\theta^2}\rho\(\frac{r-\theta}\vp\)d\theta,\ \ff r\in\rr,\ \vp>0,\label{e2.3aazi}
\end{eqnarray}
where $\rho\in C^\9_0(\rr)$, support $\rho\subset[-1,1],\ \rho\ge0$, $\int_{\rr}\rho\,dx=1$.

By \eqref{e3.4z}-\eqref{e2.3aazi},  we see that $\b_\vp\in{\rm Lip}(\rr),$ $\b'_\vp\ge0$, $b_\vp\in C^\9(\rr)$ and
\begin{eqnarray}
&\sup\{|r b_\vp(r)|;\ r\in\rr\}<\9. \label{e2.3aaazi}
\end{eqnarray}
Moreover, $b_\vp\to b$   uniformly on compacts as $\vp\to 0$.
We note that, by \eqref{e3.4z}, it follows that
\begin{equation}
\b_\vp(0)=0,\ |\b_\vp(r)|\le \inf\{|z|;\,z\in\b(r)\},\ \ \ff\vp>0,\ r\in\rr.\label{3.8a}
\end{equation}

\begin{theorem}\label{t2.2} Let $u_0\in L^1\cap L^\9.$ Then there is a unique solution $u_\vp=u_\vp(t,u_0)$ to \eqref{e2.11} which satisfies, for all $T>0$,
\begin{eqnarray}
u_\vp&\in&C([0,T];L^1)\cap C(0,T;H^1)\cap L^\9((0,T)\times\rr^d),\label{e2.12}\\
\b_\vp(u_\vp),\g_\vp (u_\vp)&\in&L^2(0,T;H^1), \label{e2.13}\\
(u_\vp)_t&\in&L^2(0,T;H^{-1}),\label{e2.14}\\
|u_\vp(t,u_0)&-&u_\vp(t,\bar u_0)|_1 \le |u_0-\bar u_0|_1,\ \ff u_0,\bar u_0\in L^1\cap L^\9\cap H^1,\qquad\ \label{e2.15}\\
|u_\vp(t)|_p&\le&|u_0|_p,\ 1\le p\le\9.\label{e2.16}
\end{eqnarray}
$u_\vp:[0,T)\to H^1$ is continuous from the right and, if $u_0\ge0$, a.e. in $\rr^d$, then $u_\vp\ge0$, a.e.  in $(0,T)\times\rr^d$. Moreover, on a subsequence $\{\vp\}\to0$, we have
\begin{eqnarray}
 u_\vp&\to& u\mbox{ in }C([0,T];L^1_{\rm loc}) ,\label{e2.16a}\\[1mm]
 \b_\vp(u_\vp)&\to&\eta^u\mbox{ weakly in }L^2(0,T;H^1),\label{e3.14ab}
\end{eqnarray}
 and weak-star in $L^\9(0,T;L^\9\cap L^2)$, where $u$ is an   entropic solution to NFPE \eqref{e1.1} with co-solution $\eta^u$. We have $u(0,x)=u_0(x)$, $\ff x\in\rr^d$ and
\begin{equation}
\label{e2.9azi}
\barr{c}
u\in C([0,T];L^1_{\rm loc}),\ \eta_u\in L^\9((0,T)\times\rr^d),\vsp
|u(t)|_p\le|u_0|_p,\ \ff p\in[1,\9],\ \ff t\in[0,T],\vsp
|u(t,u_0)-u(t,\bar u_0)|_1\le |u_0-\bar u_0|_1,\ \ff u_0,\bar u_0\in L^1,\vsp
u\ge0\mbox{, a.e. in $(0,T)\times\rr^d$, if }u_0\ge0\mbox{, a.e. in }\rr^d.\earr\end{equation}
\begin{equation}\label{e3.15ab}
\int_{\rr^d}u(t,x)dx=\int_{\rr^d}u_0(t,x)dx,\ \ff u_0\in L^1,\ u_0\ge0,\ \ff t\in[0,T].
\end{equation}
\end{theorem}

\pf To prove the  existence  for  \eqref{e2.11}, consider the operator $A_\vp$ in $H^{-1}$ defined by
$$A_\vp u={\rm div}(a_\vp(u))-\D\b_\vp(u)-\vp\D u,\ \ff u\in D(A_\vp)=H^1.$$
Arguing as in \cite{3} (see Lemma 3.1), it follows that $A_\vp$ is quasi-$m$-accretive in $H^{-1}$ and so,  there is a unique solution $u_\vp\in W^{1,\9}([0,T];H^{-1})$, which is $H^1$-valued continuous from the right on $[0,T]$.

On the other hand, the operator $A_\vp$ with the domain
 $$\left\{u\in L^1\cap H^1;\ \D\(\b_\vp(u)+ \vp u\)\in L^1\right\}$$ is accretive in $L^1$, its range contains $L^1\cap H^{-1}$, and so, its closure $\ov A_\vp$ in $L^1\times L^1$ is $m$-accretive in $L^1$.

  Here is the argument. As seen above for $f\in L^1\cap H\1$ and $\lbb$ sufficiently large, the equation
  \begin{equation}\label{2.18a}
  \lbb v_\vp+A_\vp v_\vp=f
  \end{equation}
  has a unique solution $v_\vp\in H^1$. Let $\calx_\delta\in {\rm Lip}(\rr)$, $\delta\in(0,1)$, be the following approximation of the signum function
  \begin{equation}\label{e217a}
  \calx_\delta(r)=\left\{\barr{rll}
  1&\mbox{ if }&r>\delta,\vsp
  \dd\frac r\delta&\mbox{ if }&|r|\le\delta,\vsp
  -1&\mbox{ if }&r<-\delta,\earr\right.\end{equation}
  where $\delta>0.$ The function $\calx_\delta$ is a Lipschitzian approximation of the signum function and, clearly, $\calx_\delta(v_\vp)\in H^1$. By \eqref{2.18a}, we get that
  \begin{equation}\label{2.18aa}
  \barr{r}
  \dd\lbb\int_{\rr^d}v_\vp\calx_\delta(v_\vp)dx+\int_{\rr^d}
  (\b'_\vp(v_\vp)+\vp  )
  |\nabla v_\vp|^2\calx'_\delta(v_\vp)dx\vsp
  -\dd\int_{\rr^d}(a_\vp(v_\vp)\cdot\nabla v_\vp)\calx'_\delta(v_\vp)dx
  =\int_{\rr^d}f\,\calx_\delta(v_\vp)dx.\earr\end{equation}
  On the other hand, we have
  $$\lim_{\delta\to0}\left|\int_{\rr^d}(a_\vp(v_\vp)\cdot\nabla v_\vp)\calx'_\delta(v_\vp)dx\right|
  =\lim_{\delta\to0}\left|\frac1\tau\int_{[|u_\vp|\le\delta]}a_\vp(v_\vp)\cdot\nabla v_\vp\,dx\right|=0,$$
  because $\nabla v_\vp=0$ on $\{x;\ v_\vp(x)=0\}.$ Then, letting $\delta\to0$ in \eqref{2.18aa}, we get
  $$|v_\vp|_1\le\lbb^{-1}|f|_1,\ \ff\vp>0,$$
  and so $v_\vp\in L^1$. Similarly, if $v_\vp$ and $\bar v_\vp$ are solutions to \eqref{e2.16a} corresponding to $f$ and $\bar f$ in $L^1\cap H^{-1}$, we get as above
  \begin{equation}\label{2.18aaa}
  |v_\vp-\bar v_\vp|_1\le\frac1\lbb\,|f-\bar f|_1,\ \ff\lbb>0.\end{equation}This means that the operator $A_\vp$ is accretive in $L^1$ and $L^1\cap H^{-1}\subset R(\lbb I+A_\vp)$, $\ff\lbb>0$.

  Moreover, if $f\in L^1\cap L^\9\cap H\1$, then, as easily follows by \eqref{2.18a}, we have
  \begin{eqnarray}
  |(\lbb I+A_\vp)^{-1}f|_\9&\le&\lbb^{-1}|f|_\9,\ \ff\lbb>0,\label{2.18aaaa}\\
  |(\lbb I+A_\vp)\1f|_1&\le&\lbb\1|f|_1,\ \ff\lbb>0.\label{e321aw}\end{eqnarray}
  We  note also that, if $f>0$, a.e. in $\rr^d$, then $(\lbb I+A_\vp)^{-1}f\ge0$, a.e. in $\rr^d$. (The latter follows by multiplying \eqref{2.18a} by ${\rm sign}\,v^-_\vp$  or, more exactly, by $\calx_\delta(v^-_\vp)$ and integrating over $\rr^d$.) In particular, this implies that $u_\vp\ge0$.

  Now, let us come back to equation \eqref{e2.11}.

  Multiplying by $u_\vp$, $\b(u_\vp)$   and integrating on $\delta,y)\times\rr^d$, we get
\begin{equation}\label{e2.17}
\barr{r}
\dd\frac12\,|u_\vp(t)|^2_2+\dd\int^t_0\int_{\rr^d}
\(\frac\vp2\,|\nabla u_\vp|^2
+|\nabla\b_\vp(u_\vp)|^2
\)dt\,dx\vsp
\le\dd\frac 12\,|u_0|^2_{2}+\dd\int_{\rr^d} j(u_0) dx\le C,\earr\end{equation}because
$$\barr{l}
\dd\int_{\rr^d}{\rm div}\,a_\vp(v_\vp)\b_\vp(u_\vp)dx=
\int_{\rr^d}a'_\vp(u_\vp)\b_\vp(u_\vp)\cdot\nabla u_\vp dx \dd=\int_{\rr^d}{\rm div}\,g_\vp(u_\vp)dx=0,\earr$$
and similarly for $\g_\vp(u_\vp)$.
Here $g_\vp(r)=\int^r_0a'_\vp(s)\b_\vp(s)ds,$ $j(r)=\int^r_0\b_\vp(s)ds,$  $\ff r\in\rr.$

We also note that
\begin{equation}\label{e3.21ab}
\int_{\rr^d}u_\vp(t,x)dx=\int_{\rr^d}u_0(x)dx,\ \ff t\in[0,T]    .\end{equation}
By the Crandall \& Liggett exponential formula (see \cite{1})
\begin{equation}\label{e3.23ab}
u_\vp(t)=e^{-tA_\vp}u_0=\lim_{u\to\9}\(I+\frac tn\,A_\vp\)\1u_0\mbox{ in }H\1\cap L^1,\end{equation}
it suffices to show that the solution $v\in H\1\cap L^1$ to the equation
\begin{equation}\label{e3.21aab} v+\lbb A_{\vp}v=f,\ f\in H\1\cap L^1,\ \lbb>0,\end{equation}satisfies the conservation of the mass equality
\begin{equation}\label{e3.21aaab}\int_{\rr^d}v(x)dx=\int_{\rr^d}f(x)dx.\end{equation}
By \eqref{e3.21aaab}, we have
$$\dd\int_{\rr^d} v(x)\vf(x)dx
=\dd\int_{\rr^d}f(x)\vf(x)
+\dd\int_{\rr^d}(a_\vp(v)-\nabla\b_\vp(v)-\vp\nabla v)\cdot\nabla\vf\,dx,$$
for all $\vf\in C^\9_0(\rr^d).$ Let $\mathbb{B}_N=\{x\in\rr^N;\ |x|\le N\}$ and $\vf=1$ on $\mathbb{B}_N$, $\vf=0$ on $\mathbb{B}_N^C$.
Taking into account that $a_\vp(v)\in L^1,$ $\nabla\b_\vp(v),\nabla v\in L^2$, we get, for $N\to\9$, that \eqref{e3.21aab}  holds.

Note also by \eqref{2.18aaa}--\eqref{e321aw} that it follows   that \eqref{e2.15} and \eqref{e2.16} hold.

We consider the finite difference scheme associated with equation \eqref{e2.11}, that is,
\begin{equation}\label{2.14az}
\barr{r}
\dd\frac1h\,(u^{i+1}_\vp-u^i_\vp)+{\rm div}\,a_\vp(u^{i+1}_\vp)-\D\b_\vp(u^{i+1}_\vp)
+\vp \D u^{i+1}_\vp=0,\vsp i=0,1,...,N=\mbox{$\left[\frac Th\right],$}
\earr\end{equation}
where $u^0_\vp=u_0.$ We set
\begin{equation}\label{2.19a}
u^\vp_h(t)=u^i_\vp\mbox{ for }t\in[ih,(i+1)h).
\end{equation}
By the Crandall \& Liggett formula \eqref{e3.23ab},   we know that, for $h\to0$,
\begin{equation}\label{2.19aa}
u^\vp_h\to u_\vp\mbox{ in }C([0,T];H^{-1}\cap L^1).
\end{equation}
On the other hand, since $u^0_\vp=u_0\in L^1\cap L^\9$, we have by \eqref{2.18aaa}--\eqref{2.18aaaa} and \eqref{2.19a} that $u^\vp_h\in L^\9(0,T;L^1\cap L^\9)$ and
\begin{equation}\label{2.19aaa}
|u^\vp_h|_1+|u^\vp_h(t)|_\9\le C,\ \ff t\in[0,T],
\end{equation}
where $C$ is independent of $\vp$ and $h$.

We are going to prove that, for all $\vp$ and $h$,
\begin{equation}\label{2.19aaaa}
|u^\vp_h(t+\bar t)-u^\vp_h(t)|_1\le C|\bar t|;\ t,t+\bar t\in[0,T].
\end{equation}
To this purpose, we fix $\bar t=\ell h,$ $0<\ell\le N,$ and note that, by \eqref{2.14az} and \eqref{2.19a},  we have
$$|u^{i+1+\ell}_\vp-u^{i+1}_\vp|_1\le
|u^{i+\ell}_\vp-u^i_\vp|_1\le|u^\ell_\vp-u_0|_1.$$
Then \eqref{2.19a}  follows by \eqref{2.19aaaa}, as claimed.

We note also that, by \eqref{2.19a}, it follows that
$$\int^T_0\|u^h_\vp(t)\|^2_{H^1}dt\le C,$$where $C$ is independent of $h$.

Since $H^1(\mathbb{B}_R)$ is compact in $L^1(\mathbb{B}_R)$ for all balls $\mathbb{B}_R$ of center $0$ and radius $R$, we conclude by \eqref{2.19aaaa} that along a subsequence $\{h\}\to0$ we have
$$u^h_\vp(t)\to u_\vp(t)\mbox{\ \ strongly in }L^1_{\rm loc},\ \ff\vp>0,$$
uniformly on $[0,T].$

Hence $u_\vp\in C([0,T];L^1_{\rm loc})$ and $u_\vp\in L^\9(0,T;L^1)$. By \eqref{2.19a} and the accretivity of $A_\vp$ in $L^1$, \eqref{e2.15} and \eqref{e2.16} follow.

Now, by   \eqref{e2.16a} and \eqref{e2.17}, it follows that along a subsequence converging to zero, again denoted by $\vp$, we have
\begin{eqnarray}
u_\vp&\to&u\ \ \mbox{ weak-star in }L^\9(0,T;L^2 \cap L^\9),\label{2.20}\\
\b_\vp(u_\vp)&\to&\eta\ \ \mbox{ weakly in }L^2(0,T;H^1)\label{2.21}\\
\sqrt{\vp}\,u_\vp&\to&0\ \ \mbox{ in }L^2(0,T;H^{-1}).\label{2.22}
\end{eqnarray}Moreover, since $\{u_\vp\}$ is bounded in $L^\9((0,T)\times\rr^d)$, we have also by \eqref{e3.4z}   that, for a subsequence $\{\vp\}\to0$,
\begin{equation}\label{e331a}
\eta_\vp=\b_\vp(u_\vp)  \to\eta^u\mbox{ weak-star in }L^\9((0,T)\times\rr^d).\end{equation}
Note also that, by \eqref{2.19aaaa}, we get for $h\to0$
\begin{equation}\label{2.23}
|u_\vp(t+\bar t)-u_\vp(t)|_1\le C|\bar t|,\ \ff t,t+\bar t\in[0,T],\end{equation}
and so \eqref{e2.12}  follows.
Coming back to \eqref{2.19a}, we see that, for each $\nu\in\rr^d$, we have
$$\barr{l}
\dd\frac1h\,(|u^{i+1}_\vp)_\nu-(u^i_\vp)_\nu)
+{\rm div}((a_\vp(u^{i+1}_\vp))_\nu
 -\D(\b(u_\vp^{i+1}))_\nu
)_\nu
+\vp(u^{i+1}_\vp)_\nu=0,\earr$$
where $v_\nu(x)=v(x+\nu),$ $x\in\rr^d$. Multiplying the latter by ${\rm sign}((u_\vp)^{i+1}_\nu)$ (or,~more exactly, by $\calx_\delta((u^{i+1}_\vp)_\nu))$, we get, as above,
$$|(u^{i+1}_\vp)_\nu|_1\le|(u_0)_\nu)|_1,\ \ff\nu\in\rr^d,\ i=0,1,...,N.$$Hence
\begin{equation}\label{2.24}
|(u_\vp)_\nu(t)|_1\le|(u_0)_\nu|_1,\ \ff t\in[0,T],\ \nu\in\rr^d,
\end{equation}
for all $\nu\in\rr^d$. By \eqref{2.23} and \eqref{2.24}, it follows by the Kolmogorov compactness criterium that $\{u_\vp\}_{\vp>0}$ is compact in each space $L^1(\mathbb{B}_R).$  Hence, for $\vp\to0,$
\begin{equation}\label{2.25}
u_\vp\to u\mbox{\ \ in }L^1((0,T)\times L^1_{\rm loc}).\end{equation}
In particular, it follows via Egorov's theorem that, for each $\delta>0$, there is a measurable subset $\Sigma_\delta\subset(0,T)\times\mathbb{B}_R$ whose complement in $(0,T)\times\mathbb{B}_R$ has   Lebesgue measure less than $\delta$ such that $u_\vp\to u$ in $L^\9(\Sigma_\delta)$. Since  $\b:L^2(\Sigma_\delta)\to L^2(\Sigma_\delta)$ is maximal monotone, strongly-weakly closed, and so  it follows by \eqref{2.21} and \eqref{e331a} that $\eta^u\in\b(u)$ a.e. in $\Sigma_\delta$ and, since $\delta$ is arbitrary, it follows   that $\eta^u\in\b(u)$, a.e. in $(0,T)\times\rr^d.$ Moreover, by \eqref{2.21}, it follows  that $\nabla \b_\vp (u_\vp)\to\nabla\eta^u$ weakly in $L^2(0,T;L^2)$, as $\vp\to0$.

Note also that, letting $\vp\to0$ in \eqref{e2.15}--\eqref{e2.16}, it follows by \eqref{2.25} via Fatou's lemma that \eqref{e2.9azi} holds.

Since by   \eqref{2.23}, $\{u_\vp\}$ is equi-uniformly  continuous in $C([0,T];L^1)$, by  \eqref{2.25} it follows that $u_\vp\to u$ in $C([0,T];L^1_{\rm loc})$ as $\vp\to0$. In particular, this implies that $u\in C([0,T];L^1_{\rm loc})$. Letting $\vp\to0$ in \eqref{e3.21ab}, we get also \eqref{e3.15ab}.

Let us prove now that $u$ is an entropic solution for equation \eqref{e1.1}, that is, \eqref{e2.1prim} holds.
By \eqref{e2.11}, we have
\begin{equation}\label{2.26}
 \barr{l}
\dd\int^T_0\!\!\!\int_{\rr^d}\calx_\delta(u_\vp-k)
\!\left[(u_\vp)_t\vf+{\rm div}(a_\vp(u_\vp))\vf
+\vp u_\vp\vf\right]dx\,dt\\
\qquad\quad+\dd\int^T_0\!\!\!\int_{\rr^d}
(\nabla\b_\vp(u_\vp))\cdot\nabla (\calx_\delta(u_\vp-k)\vf)dx\,dt=0.
\earr\!\!\!\end{equation}Now, letting $\vp\to0$ in \eqref{e3.21ab}, we get \eqref{e3.15ab}, as claimed.

We have, for $j_\delta(r)\equiv\int^r_0\calx_\delta(s)ds,$
\begin{equation}\label{2.27}
 \barr{l}
\dd\int^T_0\int_{\rr^d}\calx_\delta(u_\vp-k)(u_\vp)_t\vf\,dx\,dt
\vsp
\qquad=-\dd\int^T_0\int_{\rr^d}\calx_\delta(u_\vp-k)(u_\vp-k)_t\vf\,dx\,dt\vsp
\qquad=\dd\int^T_0\int_{\rr^d}j_\delta(u_\vp-k)
\vf_tdx\,dt
\buildrel{\delta\to0}\over{\longrightarrow}
\dd\int^T_0\int_{\rr^d}|u_\vp-k|\vf_tdx\,dt.\earr
\end{equation}
We also have
\begin{equation}\label{2.28}
\barr{r}
\left|\dd\int^T_0\!\!\!\int_{\rr^d}
\calx_\delta(u_\vp-k)
(\vp u_\vp\vf)dx\,dt\right|
\le C\sqrt{\vp}\,\|\vf\|_{L^1((0,T)\times\rr^d)},
\earr\end{equation}
and
\begin{equation}\label{2.29}
\barr{l}
\dd\int^T_0\!\!\!\int_{\rr^d}
\nabla\b_\vp(u_\vp)\cdot\nabla(\calx_\delta(u_\vp-k)\vf)dx\,dt\vsp
\qquad
=\dd\int^T_0\!\!\!\int_{\rr^d}
\b'_\vp(u_\vp)|\nabla u_\vp|^2\calx'_\delta(u_\vp-k)\vf\,dx\,dt\vsp
\qquad
+\dd\int^T_0\!\!\!\int_{\rr^d}
(\nabla\b_\vp(u_\vp)\cdot\nabla\vf)\calx_\delta(u_\vp-k)dx\,dt\vsp
\qquad
\ge\dd\int^T_0\!\!\!\int_{\rr^d}
(\nabla\b_\vp(u_\vp)\cdot\nabla\vf)\calx_\delta(u_\vp-k)dx\,dt,
\earr\end{equation}
because $\vf\ge0$. Furthermore,
\begin{equation}\label{2.30}
\barr{l}
\dd\int^T_0\!\!\!\int_{\rr^d}
\calx_\delta(u_\vp-k){\rm div}(a_\vp(u_\vp))\vf\,dx\,dt\vsp
\qquad=
\dd\int^T_0\!\!\!\int_{\rr^d}
\calx_\delta(u_\vp-k){\rm div}(a_\vp(u_\vp)-a_\vp(k))\vf\,dx\,dt\vsp
\qquad=
\dd\int^T_0\!\!\!\int_{\rr^d}
\calx_\delta(u_\vp-k){\rm div}(a_\vp(u_\vp)-a_\vp(k))\vf\,dx\,dt\vsp
\qquad+
\dd\int^T_0\!\!\!\int_{\rr^d}
\calx_\delta(u_\vp-k)(a_\vp(u_\vp)-a_\vp(k))\cdot\nabla\vf\,dx\,dt
=I^1_{\vp,\delta}+I^2_{\vp,\delta}.\earr
\end{equation}
We note that
$$\barr{ll}
\calx_\delta(u_\vp-k){\rm div}(a_\vp(u_\vp)-a_\vp(k))
\!\!\!&={\rm div}(\calx_\delta(u_\vp-k)(a_\vp(u_\vp)-a_\vp(k))\vsp
&-\dd\int^{u_\vp}_k\calx'_\delta(s)(a_\vp(u_\vp(s))-a_\vp(k))ds).\earr$$
This yields
$$\barr{lcl}
I^1_{\vp,\delta}&=&-
\dd\int^T_0\!\!\!\int_{\rr^d}\nabla\vf
\cdot(\calx_\delta(u_\vp-k)
(a_\vp(u_\vp)-a_\vp(k))dx\,dt
\vsp
&&-\dd\int^{u_\vp}_k\calx'_\delta(s)(a_\vp(u_\vp(s)-a_\vp(k))ds)dx\,dt\earr$$
and so, we have
\begin{equation}\label{2.31}
\lim_{\vp,\delta\to0}I^1_{\vp,\delta}=0.\end{equation}
We have also, by Lebesgue's dominated convergence theorem, that
\begin{equation}\label{2.32}
\lim_{\vp,\delta\to0}I^2_{\vp,\delta}=
\dd\int^T_0\!\!\!\int_{\rr^d}{\rm sign}(u-k)(a(u)-a(k))\cdot\nabla\vf\,dx\,dt.\end{equation}
Now, combining \eqref{2.26}--\eqref{2.32} and letting $\vp,\delta\to0$, we obtain that
$$\barr{ll}
-\dd\int^T_0\!\!\!\int_{\rr^d}(u-k)\vf_tdx\,dt\!\!\!
&-\dd\int^T_0\!\!\!\int_{\rr^d}
{\rm sign}(u-k)(a(u)-a(k))\cdot\nabla\vf\,dx\,dt\vsp
&+\dd\int^T_0\!\!\!\int_{\rr^d}
\nabla\eta^u\cdot\nabla\vf\,{\rm sign}(u-k)dx\,dt\le0,\earr$$
and so \eqref{e2.1prim} follows.

  Then it follows that $u$ is an entropic solution to \eqref{e1.1} with co-solution $\eta^u$. This completes the proof.

\begin{remark}\label{r2.1}{\rm It should be mentioned that, if $a$ is not Lipschitz, then the entropic solution $u$ to \eqref{e1.1} in general is not unique (this happens also for the special case of the   conservation law equation, that is, for $\b\equiv0$). On the other hand, the solution $u=\lim\limits_{\vp\to0}u_\vp$ given by Theorem \ref{t2.2} is just the mild solution (see \cite{1}, p. 128) to the Cauchy problem
\begin{equation}\label{2.33}
\barr{l}
\dd\frac{du}{dt}+Au\ni0,\ t\in(0,T),\vsp
u(0)=u_0,\earr\end{equation}
where $A$ is the closure in $L^1\times L^1$ of the operator $A_0:D(A_0)\subset L^1\to L^1$ defined by $v\in A_0u$  if $u,v\in L^1$, $a(u)\in L^1$, $\eta\in H^1_{\rm loc}$,
$\eta\in\b(u)$, a.e. in $\rr^d$, and
$$\int_{\rr^d}|{\rm sign}(u-k)((a(u)-a(k))\cdot\nabla\vf
-\nabla(\eta\cdot\nabla\vf)+v\vf)dx\ge0,$$
for all $k\in\rr$ and $\vf\in C^\9_0(\rr^d),$ $\vf\ge0$.}\end{remark}

It turns out (see \cite{2}) that $A$ is $m$-accretive in $L^1$ and so, by the Crandall \& Liggett theorem, \eqref{2.33} has a unique mild solution $y=e^{-tA}y_0=\lim\limits_{n\to\9}\(I+\frac tn\,A\)^{-n}y_0$ in $L^1.$
On the other hand, arguing as in the proof of Theorem \ref{t2.2}, it follows that, for each $f\in L^1\cap L^\9$, the solution $u^\vp$ to the equation $\lbb u^\vp+A_\vp u^\vp=f$ is strongly convergent in $L^1_{\rm loc}$ for $\vp\to0$ and $\lbb>0$ to $v=(\lbb I+A_0)^{-1}f$. This implies by an argument similar to that used in the proof of the Trotter--Kato theorem for nonlinear semigroups (see, e.g., \cite{1}, p.169) that, for $u_0\in L^1\cap L^\9$, the solutions $u_\vp$ to equation \eqref{e2.11} are convergent in $C([0,T];L^1)$ to $y=e^{-tA}y_0$.  Hence, we can identify in the class of entropic solutions to \eqref{e1.1} the function $u$ given by Theorem \ref{t2.2}  as the evolution generated by the operator $A$ in $L^1$. Such a solution will be called in the following the {\it generalized entropic solution} to FPNE \eqref{e1.1}.

\section{The  probabilistic representation\\ of generalized entropic solution}
\setcounter{equation}{0}
\setcounter{theorem}{0}

We set
\begin{equation}\label{e3.2}
\Phi(r)=\(\frac{\b(r)}r\,\one_{[0,\9)}(r)\)^{\frac12},\ \ff r\in\rr,\end{equation}
and assume that, besides (i),    the following conditions hold
\bit\item[(iii)] $\sup\{z;\ z\in\Phi(r)\}\le\rho<+\9,\ \ff r\in[0,\9)$. \eit

Consider the Wiener process
\begin{equation}\label{e3.3}
W={\rm col}(W_j)^d_{j=1},\end{equation}
where $W_j$ are independent Brownian motions on a filtered probability space $(\ooo,\calf,(\calf_t),\mathbb{P}).$
We note that, by virtue of assumptions (ii) and (iii), the function $\Phi:\rr\to\rr$ is continuous and bounded.

Let $u=u(t,x)$ be the generalized entropic solution to \eqref{e1.1} with the corresponding co-solution $\eta^u$. We set
$$\calx^u(t,x)=\left\{\barr{ll}
\(\dd\frac{\eta^u(t,x)}{u(t,x)}\)^{\frac12}&\mbox{ on }[(t,x);u(t,x)>0],\vsp
0&\mbox{ on }[(t,x);u(t,x)=0],\earr\right.$$
and note that $\calx^u$ is measurable and bounded on $(0,T)\times\rr^d$.

Consider the stochastic differential equation
\begin{equation}\label{e3.4}
\barr{l}
dY(t)=b(u(t,Y(t)))dt+\sqrt{2}\,\calx^u(t,u(t,Y(t)))
 dW(t),\ t\ge0,\vsp
Y(0)=Y_0.\earr\end{equation}

A probabilistically weak solution $Y : [0,T]\times\Omega \to \rr^{+}$ to
equation \eqref{e3.4} is a  progressively measurable and  pathwise continuous process which satisfies \eqref{e3.4} in integral form for some $\rr^d$-valued Wiener process on some filtered probability space $(\ooo,\calf,(\calf)_t,P)$ as above. Theorem \ref{t3.1} below is the main result of this work.

\begin{theorem}\label{t3.1} Let assumptions {\rm(i)-(iii)} hold and let $u$ be the generalized entropic solution to NFPE \eqref{e1.1} given by Theorem~{\rm\ref{t2.2}}  with $u_0\ge0$ in $\rr^d$, $u_0\in L^1\cap L^\9$, $\int_{\rr^d}u_0dx=1
	$.  Then there is a probabilistically weak solution $Y$   to \eqref{e3.4}   such that
\begin{equation}\label{e3.5}
P\circ Y(t)^{-1}(dx)=u(t,x)dx\mbox{ for every }t\in[0,T].  \end{equation}
\end{theorem}

\pf We recall that $u$ is also a distributional solution and also that the function $t\to u(t,x)$ is $L^1_{\rm loc}$ continuous on $[0,T]$, $u\in L^1((0,T)\times\rr^d)\cap L^\9((0,T)\times\rr^d)$ and $\int_{\rr^d}u(t,x)dx=1$, $\ff t\in[0,T]$.  Moreover, by (iii) we have also
$$\int^T_0\int_{\rr^d}(|b(u(t,x))|+
|\eta^u(t,x)|)
|u(t,x)|dt\,dx<\9.$$
Then, as explained in Section 2, by \cite{15},  Theorem 2.5, and \cite{14prim}, Theorem 2.6, there exists a $d$-dimensional $(\calf_t)_{t\ge0}$-Brownian motion $W(H)$, $t\in[0,T]$, on a stochastic basis $(\ooo,\calf,(\calf_t)_{t\ge0},\mathbb{P})$ and a continuous   $(\calf_t)$-progressively measurable map $Y:[0,T]\times\ooo\to\rr^d$ which satisfies equation \eqref{e3.4},
and $P_{u_0}\circ Y(t)^{-1}(dx)=u(t,x)dx$ for all \mbox{$t\in[0,T].$}
This completes the proof.

\section{NFPE with linear drift; distributional\\ solutions}
\setcounter{equation}{0}
\setcounter{theorem}{0}

Consider the equation
\begin{equation}\label{e5.1}
\barr{l}
u_t+{\rm div}(D(x)u)-\D\b(u)=0\mbox{ in }(0,\9)\times\rr^d,\vsp
u(0,x)=u_0(x),\earr\end{equation}
where $D:\rr^d\to\rr^d$   satisfies the following condition
\bit\item[(k)] $D\in L^\9(\rr^d),\ {\rm div}\,D\in L^\9(\rr^d)$.
\eit
It turns out (see, e.g., \cite{2a}) that, under assumptions (i) and (k), for $u_0\in L^1\cap L^\9$, equation \eqref{e5.1} has a unique solution $u\in C([0,T];L^1)\cap L^\9((0,T)\times\rr^n)$ in sense of distributions, that is,
\begin{equation}\label{e5.2}
\barr{l}
u_t+{\rm div}(Du)-\D\b(u)=0\mbox{ in }\cald'((0,T)\times\rr^d),\vsp
u(0)=u_0.\earr\end{equation}
Moreover, $u\ge0$ if $u_0\ge0$ and $\int_{\rr^d}u(t,x)dx=\int_{\rr^d}u_0(x)dx,\ \ff t\in[0,T].$

The solution $u$ is obtained as in the previous case as limit of solutions $u_\vp$ to the approximating equation
$$\barr{l}
(u_\vp)_t+{\rm div}(Du_\vp)-\D(\b_\vp(u_\vp)+\vp u_\vp)=0\mbox{ in }(0,T)\times\rr^d,\vsp
u_\vp(0)=u_0.\earr$$

\begin{theorem}\label{t5.1} Let $u$ be the solution to \eqref{e5.2} for $u_0\in L^1\cap L^\9\cap H^1\cap W^{1,\9}_{\rm loc},$ $u_0\ge0$. Then
\begin{equation}\label{e5.3}
u=\call(Y),\mbox{ a.e. in }(0,T)\times\rr^d,\end{equation}
where $\call(Y)$ is the probability density law of a weak solution $Y$ to the stochastic differential equation
\begin{equation}\label{e4.3a}
\barr{l}
dY= D(Y(t))dt+\calx^u(t,Y(t))dW,\vsp
Y(0)=Y_0,\earr\end{equation}
where $\calx^u$ is as in Section {\rm4}.
\end{theorem}

The proof follows as in the case of Theorem \ref{t3.1} by the results of Section~2.

\begin{remark}\rm \label{r4.1} We note also that the bosons equations   $(\b (u)=\g\ln(u+1))$   enters within the applicability field of Theorem \ref{t5.1}. As a matter of fact, it applies to each monotone continuous function $\b$ with polynomial growth, $\b(r)\le Cr^m$, $\ff r\in\rr$. Also, the self-organized criticality model presented in Section 1 is covered by Theorems \ref{t3.1} and \ref{t5.1}. \end{remark}

\section{An application to $1{-}D$ degenerate parabolic equation}
\setcounter{equation}{0}
\setcounter{theorem}{0}

Consider the equation
\begin{equation}
\label{e6.1}\barr{l}
y_t(t,x)+a(y_x(t,x))-(\b(y_x(t,x)))_x=0,\ t\ge0,\ x\in\rr,\vsp
y(0,x)=y_0(x),\ x\in\rr,\earr
\end{equation}where $a(r)\equiv b(r)r$ and $b:\rr\to\rr,$ $\b:\rr\to\rr$ satisfy assumptions (i), (ii), (iii). By the substitution $u=y_x$, we reduce equation \eqref{e6.1} to a Fokker-Planck equation of the form \eqref{e1.1}, namely
\begin{equation}
\label{e6.2}
\barr{l}
u_t+(a(u))_x-(\b(u))_{xx}=0,\ t\ge0,\ x\in\rr,\vsp
u(0,x)=y'_0(x),\ x\in\rr.\earr\end{equation}By Theorem \ref{t2.2}, equation \eqref{e6.2}  has an entropic solution $u$ which sa\-tis\-fies \eqref{e2.9azi}.  Then
$$y(t,x)=\int^x_{-\9}u(t,\xi)d\xi,\ t\ge0,\ x\in\rr,$$can be viewed as a generalized solution to \eqref{e6.1}. Moreover, by Theorem \ref{t3.1}, it follows that, if
\begin{equation}
\label{e6.3}
y'_0\ge0,\ y_0(+\9)-y_0(-\9)=1,
\end{equation}then $y$  has the probabilistic representation
\begin{equation}
\label{e6.4}
y(t,x)=y(t,0)+\mathbb{P}[0\le Y(t)\le x],\ \ff t\ge0,\ x\in\rr,
\end{equation}where $Y$ is the solution to the stochastic differential equation \eqref{e3.4}.

\bk\n{\bf An example.} The equation
\begin{equation}
\label{e6.5}
\barr{l}
y_t+a(y_x)-a_1(y_x)(U'(y_x))_x=0,\ t>0,\ x\in\rr,\vsp
y(0,x)=y_0(x),\ x\in\rr,\earr
\end{equation}where $U:\rr\to\rr$ is a convex, and smooth function, $U(0)=0$, $a_1\in C_b(\rr),$ $a_1\ge0$ and $a\in C(\rr)$, was proposed as a model for the evolution law of interfacial curves arising in material science and crystal  growth. (See \cite{11z}) for its physical significance and an analytic treatment via comparison arguments.)

By the substitution
$$\b(r)=\int^r_0 a_1(s)U''(s)ds,\ \ff r\in\rr,$$ one reduces \eqref{e6.5} to equation \eqref{e6.1} and so, if $U''\in L^\9(\rr)$, we have for the solution $y$ to \eqref{e6.5} the probabilistic representation formula \eqref{e6.4}, where $Y$ is he solution to SDE
$$dY=\frac{a(Y)}Y\,dt+\(\frac2Y\,\int^Y_0a_1(s)U''(s)\)dW.$$
As far as concerns the existence for \eqref{e6.5}, the condition $U'\in W^{1,\9}(\rr)$ can be weakened to $U\in W^{1,\9}(\rr)$, which leads to a maximal monotone graph $\b:\rr\to2^\rr$ everywhere defined on $\rr$.

\bk\n{\bf Acknowledgement.} This work was supported by the DFG through CRC
1283. V.~Barbu was partially supported by the CNCS-UEFISCDI project
PN-III-P4-ID-PCE 2015-0011.


\begin{thebibliography}{nn}

\bibitem{1} V. Barbu, {\it Nonlinear Differential Equations of Monotone Type in Banach Spaces}, Springer, 2010.

\bibitem{2} V. Barbu, Generalized solutions to nonlinear Fokker-Planck equations, {\it J. Differential Equations}, 261 (2016), 2446-2471.

\bibitem{2a} V. Barbu, Generalized solutions to nonlinear Fokker-Planck
equations with linear drift, {\it Stochastic Partial Differential Equations and Related Fields -- In Honor
of Michael R\"ockner}, Springer Proceedings in Ma\-the\-ma\-tics and
Statistics, 2017.


\bibitem{3ab} V. Barbu, The steepest descent algorithm in Wasserstein metric for the sandpile model of self-organized criticality, {\it SIAM J. Control Optim.}, 56 (2017), 413-422.

\bibitem{3} V. Barbu, M. R\"ockner, Nonlinear Fokker-Planck equations driven by Gaussian linear multiplicative noise, arXiv:1708.08768.

\bibitem{4} V. Barbu, M. R\"ockner, F. Russo, Probabilistic representation for solutions of an irregular porous media type equation; the degenerate case, {\it Probab. Theory Rel. Fields}, 15 (2011), 1-43.

\bibitem{5}   Ph. Blanchard, M. R\"ockner, F.Russo, Probabilistic representation forf solutions of an irregular porous media type equation, {\it The Annals of Probability}, 38 (2010), 1870-1900.

    \bibitem{7a} S. Benachour, P. Chassaing, B. Roynette, P. Vallois, Processus associ\'ees \`a l'\'equation de millieux poreux, {\it Ann. Scuola Norm. Sup. Pisa Cl.  Sci}.,   23 (4) (1996), 793-838.


\bibitem{6}   M.G. Crandall, The semigroup approach to first order quasilinear equations in several space variables, {\it Israel J. Math.}, 12 (1972), 108-132.

\bibitem{9a} A. Figalli, Existence and uniqueness of martingale solutions for SDEs with rough or degenerate coefficients, {\it J. Funct. Anal.}, 254 (2008), 109-153.

\bibitem{7} T.D. Frank, {\it Nonlinear Fokker-Planck Equations}, Springer-Verlag, 2005.


\bibitem{8} T.D. Frank, A. Daffertshofer, $H$-Theorem for nonlinear Fokker--Planck equations related to generalized thermostatistics, {\it Physica A}, 295 (2001), 455-477.


\bibitem{11z} M.H. Giga, Y. Giga, Evolving graphs by singular weighted curvature, {\it Archive Rational Mech. Anal.}, 141 (1998), 117-198.

\bibitem{9} S.N. Kru\v zkov, First order quasilinear equations in several independent variables, {\it Math. USSR-Sbornik}, 10 (1970), 217-243.

\bibitem{13a} A. Le Cavil, N. Oudjane, F. Russo, Probabilistic representation of a class of non conservative nonlinear partial differential equations, {\it ALEA Lat. Amer. J. Probab. Math. Statistic}, 13 (2016), 1189-1233.


\bibitem{15a} M. R\"ockner, F. Russo, Uniqueness for a class of stochastic Fokker-Planck and porous media equations, {\it J. Evol. Equ.}, 17 (3) (2017), 1049-1062.

\bibitem{14prim} D.W. Stroock, {\it Lectures on stochastic analysis: difusion theory}, LMS Student Texts, 1987.




\bibitem{14} D.W. Stroock, S.R.S. Varadhan, {\it Multidimensional diffusion processes}, Springer, 1979.




\bibitem{15} D. Trevisan, Well-posedness of multidimensional diffusion processes with weakly differentiable coefficients, {\it Electron J. Probab.}, 21 (22) (2016), 1-41.

\bibitem{16} F.Y. Wang, Distribution-dependent SDEs for Landau type equations, arXiv:1606.05843. To appear in {\it Stoch. Proc. Appl.}



\end{thebibliography}
\end{document}